\documentclass[final,twoside,11pt]{entics}
\usepackage{enticsmacro}
\usepackage{graphicx}
\usepackage[all]{xy}
\newcommand{\ua}{\mathord{\uparrow}}
\newcommand{\da}{\mathord{\downarrow}}
\newcommand{\rom}[1]{\rm{\uppercase\expandafter{\romannumeral #1}}}

\def\conf{ISDT 2022} 	
\volume{2}			
\def\volu{2}			
\def\papno{14}			
\def\mydoi{{\normalshape  \href{https://doi.org/10.46298/entics.10369}{10.46298/entics.10369}}}
\def\copynames{L.. Zhang X. Zhou, Q. Li} 
\def\runtitle{A Hofmann-Mislove theorem for $c$-well-filtered spaces}

\def\lastname{Zhang, Zhou and Li}
\begin{document}
\begin{frontmatter}
  \title{A Hofmann-Mislove Theorem for $c$-well-filtered Spaces}
  \author{Liping Zhang\thanksref{myemail}}	
   \author{Xiangnan Zhou\thanksref{ALL}\thanksref{coemail}}
   \author{Qingguo Li\thanksref{ALL}\thanksref{coemaill}}		
   \address{School of Mathematics\\Hunan University\\ 
             Changsha, Hunan, 410082, China}  							
  \thanks[ALL]{This work is supported by the National Natural Science Foundation of China (Grant No. 12231007) and the Natural Science Foundation of Hunan Province (Grant No. 2019JJ50041).}   
   \thanks[myemail]{Email: \href{mailto:myuserid@mydept.myinst.myedu} {\texttt{\normalshape
        lipingzhang@hnu.edu.cn}}}
  \thanks[coemail]{Corresponding author. Email:  \href{mailto:couserid@codept.coinst.coedu} {\texttt{\normalshape
        xnzhou81026@163.com, xnzhou@hnu.edu.cn}}}
  \thanks[coemaill]{Email:  \href{mailto:couserid@codept.coinst.coedu} {\texttt{\normalshape
        liqingguoli@aliyun.com}}}
\begin{abstract}
 The Hofmann-Mislove theorem states that in a sober space, the nonempty Scott open filters of its open set lattice  correspond bijectively to its compacts saturated sets. In this paper, the concept of $c$-well-filtered spaces is introduced. We show that a retract of a $c$-well-filtered space is $c$-well-filtered and a locally Lindel\"{o}f and $c$-well-filtered $P$-space is countably sober. In particular, we obtain a Hofmann-Mislove theorem for $c$-well-filtered spaces.
\end{abstract}
\end{frontmatter}
\section{Introduction}\label{intro}
The Hofmann-Mislove Theorem plays an important role in the study of the basic topological theorems concerning sober spaces and illustrates the close relationship between domain theory and topology. It states that there exists a bijection between the nonempty Scott open filters on the open set lattice for a sober space $X$ and the compact saturated subsets (\cite{Gierz03,Hofmann81}). Moreover, it also declares that there is a bijection between the family of nonempty Scott open filters of the compact saturated sets and the open set lattice in a locally compact sober space $X$.

In recent years, some researchers have generalized the Hofmann-Mislove Theorem to some other topological spaces (\cite{jung2007,kovar2005,schroder2015,xu2009,yu2021,yang17}). For example, A.~Jung gave an analogy result of the Hofmann-Mislove theorem for bisober spaces in~\cite{jung2007}. In~\cite{schroder2015}, M.~Schr\"{o}der proved that there exists a continuous retraction from the family $\mathcal {O}\mathcal {O}_{+}(X)$ of all nonempty $\sigma$-Scott-open collections of open sets to the upper space $\mathcal {K}(X)$ of all countably-compact sets in a sequentially Hausdorff sequential space $X$. J.B.~Yang and J.M.~Shi obtained that in a countably sober space, a Scott open countable filter of open set lattice is precisely a compact saturated set in~\cite{yang17}. The motivation of this paper is to establish the relationship between the open set lattice and the set of all $\sigma$-Scott-open countable filters of saturated Lindel\"{o}f sets in a $c$-well-filtered space.

The remaining parts of this paper are organized as follows. Section 2 recalls some basic concepts and results used in this paper. In Section 3, we define a new notion of $c$-well-filtered spaces and investigate its basic properties. Particularly, we prove that $c$-well-filteredness is hereditary for saturated subsets and a retract of a $c$-well-filtered space is $c$-well-filtered. In Section 4, we obtain a Hofmann-Mislove Theorem for $c$-well-filtered spaces, which states that there is a bijection between the open set lattice and the set of all $\sigma$-Scott-open countable filters of saturated Lindel\"{o}f sets.

\section{Preliminaries}
We refer to~\cite{Gierz03,Han89,Lee88} for the standard definitions and notations of order theory and domain theory, and to~\cite{goubault13,Gillman76,McGovern01} for topology.

Let $X$ be a set. We denote the family of all finite subsets (resp., countable subsets) of $X$ by Fin$X$ (resp., Count$X$). Let $L$ be a poset. A nonempty subset $D\subseteq L$ is \emph{countably directed} if for every $E\in\mbox{Count}D$, there exists $d\in D$ such that $E\subseteq\da d$. \emph{Countably filtered} is defined dually. A nonempty subset $F\subseteq L$ is called a \emph{countable filter} if it is a countably filtered upper set. A poset $L$ is called a \emph{countably directed complete poset} if every countably directed subset $D\subseteq L$ has a least upper bound $\sup D$ in $L$. An upper set $U$ of $L$ is \emph{$\sigma$-Scott open} if for every countably directed subset $D\subseteq L$, $\sup D\in U$ implies $D\cap U\neq\emptyset$. All $\sigma$-Scott open subsets of $L$ form a topology, called the \emph{$\sigma$-Scott topology} and denoted as $\sigma_{c}(L)$.

Let $L$ be a poset. The symbol $\sigma(L)$ denotes the \emph{Scott topology} consisting of all Scott open subsets of $L$. The space $\Sigma L=(L, \sigma(L))$ is called the \emph{Scott space} of $L$.


Let  $X$ be a topological space and $\mathcal {O}(X)$ the open set lattice. A subset $A$ of $X$ is a \emph{Lindel\"{o}f} set if each open cover $\mathcal {U}$ of $A$ has a countable subcover. We denote the set of all compact saturated (resp., saturated Lindel\"{o}f) subsets of $X$ by $\mathcal {Q}(X)$ (resp., $\mathcal {L}\mathcal {Q}(X)$). A topological space $X$ is \emph{well-filtered} iff for every filtered family $\mathcal {K}$ of $\mathcal {Q}(X)$ and for every open subset $U$ of $X$, if $\bigcap\mathcal {K}\subseteq U$, then $K\subseteq U$ for some $K\in\mathcal {K}$. A topological space $X$ is \emph{locally Lindel\"{o}f} if for every open subset $U$ of $X$ and for every point $x\in U$, there exists $K\in\mathcal {L}\mathcal {Q}(X)$ such that $x\in \mbox{int}(K)\subseteq K\subseteq U$. A point $p\in X$ is called a \emph{$P$-point} if its neighbourhood system is closed under countable intersection. A topological space $X$ is called a \emph{$P$-space} if every point in $X$ is a $P$-point.

In what follows, we will give some results on Lindel\"{o}f sets.

\begin{proposition}\label{prop-2.1}
Let $X$ be a topological space. A subset $K$ of $X$ is a Lindel\"{o}f set if and only if for every countably directed family $\{U_{i}\}_{i\in I}$ of open subsets of $X$, if  $K\subseteq\bigcup_{i\in I}U_{i}$, then $K\subseteq U_{i}$ for some $i\in I$.
\end{proposition}
\begin{proof}
(If part) Let $\mathcal {U}$ be an open cover of $K$. Take $\mathcal {V}=\{U_{\alpha}: U_{\alpha}=\bigcup_{i\in \mathbb{N}}U_i, U_i\in \mathcal{U}, \alpha\in I\}$, where $I $ is an index set. Then $\mathcal {V}$ is a countably directed family of open subsets of $X$ and $\bigcup\mathcal {V}=\bigcup\mathcal{U}$. Thus there exists $\alpha\in I$ such that $K\subseteq U_{\alpha}$, which means $\mathcal {U}$ has a countable subcover.

(Only if part) It is obvious by the countably directedness of $\{U_{i}\}_{i\in I}$.
\end{proof}

\begin{proposition}\label{prop-2.2}
Let $X$ be a topological space. A subset $K$ of $X$ is a Lindel\"{o}f set if and only if for every countably filtered family $\mathcal {C}$ of closed subsets of $X$, if $K\cap C\neq\emptyset$ for all $C\in\mathcal {C}$, then $K\cap\bigcap\mathcal {C}\neq\emptyset$.
\end{proposition}

By Proposition \ref{prop-2.2}, we get immediately the following corollary.

\begin{corollary}\label{coro-2.1}
Let $X$ be a topological space. If $K$ is a Lindel\"{o}f subset of $X$ and $C$ is a closed subset of $X$, then $K\cap C$ is a Lindel\"{o}f set.
\end{corollary}

\begin{proposition}\label{prop-2.3}
Let $f: X\rightarrow Y$ be a continuous map. If $C$ is a Lindel\"{o}f subset of $X$, then $f(C)$ is a Lindel\"{o}f set of $Y$.
\end{proposition}

\begin{proposition}\label{prop-2.4}
Let $X$ be a Lindel\"{o}f space. Then every closed subset $F$ of $X$ is a Lindel\"{o}f set.
\end{proposition}

\section{$c$-well-filtered spaces}
In this section, we introduce a notion of $c$-well-filtered spaces and discuss its some basic properties.

\begin{definition}
A topological space $X$ is \emph{$c$-well-filtered} if for every countably filtered family $\{K_{i}\}_{i\in I}$ of saturated Lindel\"{o}f subsets of $X$ and each open subset $U$ with $\bigcap_{i\in I}K_{i}\subseteq U$, there is a $K_{i_{0}}\subseteq U$ for some $i_{0}\in I$ .
\end{definition}

\begin{proposition}\label{prop-3.1}
Suppose that $X$ is a countable set. Then every topological space $(X, \tau)$ is $c$-well-filtered.
\end{proposition}
\begin{proof}
Obviously, each subset of $X$ is a Lindel\"{o}f set. Hence $\mathcal {L}\mathcal {Q}(X)=\{{\ua}K: K\subseteq X\}$. Let $\mathcal {K}$ be a countably filtered family of $\mathcal {L}\mathcal {Q}(X)$ and $U\in \tau$ with $\bigcap\mathcal {K}\subseteq U$. For every nonempty chain $\mathcal {C}$ in $\mathcal {K}$, there exists $K\in \mathcal {K}$ such that $K$ is a lower bound of $\mathcal {C}$ by the countably filteredness of  $\mathcal {K}$ .
By the order-dual of Zorn's Lemma, $\mathcal {K}$ contains a minimal element $K_{0}$. Hence, $K_{0}=\bigcap\mathcal {K}\subseteq U$. Therefore, $(X, \tau)$ is $c$-well-filtered.
\end{proof}

The following example shows that $c$-well-filtered spaces are not always well-filtered spaces.

\begin{example}\label{exam-1}
Let $\mathbb{J}=\mathbb{N}\times(\mathbb{N}\cup\{\omega\})$ with the partial order defined by $(j, k)\leq(m, n)$ iff $j=m$ and $k\leq n$, or $n=\omega$ and $k\leq\omega$, where $\mathbb{J}$ is a well-known dcpo constructed by Johnstone in \cite{Johnstone81}. It is easy to see that $\mathbb{J}$ is countable. By Proposition~\ref{prop-3.1},  $\Sigma\mathbb{J}$ is a $c$-well-filtered space. But $\Sigma\mathbb{J}$ is not a well-filtered space (see~\cite[Example 3.1]{lu17}).
\end{example}

Conversely, a well-filtered space may not be a $c$-well-filtered space.

\begin{example}\label{exam-2}
Consider the real number set $\mathbb{R}$ with the co-countable topology $\tau_{coc}$, where $\tau_{coc}=\{U\subseteq\mathbb{R} :\mathbb{R}\backslash U~\mbox{is countable}\}\bigcup \{\emptyset\}$. It is known that  the topological space $(\mathbb{R}, \tau_{coc})$ is a well-filtered $T_{1}$-space (see~\cite[Example 3.14]{shan22}).

Next, we show that all subsets of $\mathbb{R}$ are saturated Lindel\"{o}f sets.
Let $K$ be a subset of $\mathbb{R}$ and assume that $\mathcal {U}$ is an open cover of $K$. Given $U\in\mathcal {U}$, then there exists a countable subset $C$ of $\mathbb{R}$ such that $U=\mathbb{R}\backslash C$ and $K=(K\cap C)\cup(K\backslash C)\subseteq\bigcup\mathcal {U}$. It is obvious that $K\backslash C\subseteq\mathbb{R}\backslash C=U$. Since $K\cap C$ is countable and $K\cap C\subseteq\bigcup\mathcal {U}$, there exist countably many members of $\mathcal {U}$ whose union contains $K\cap C$. Hence, $K$ is a saturated Lindel\"{o}f set.

Finally, we show that $(\mathbb{R}, \tau_{coc})$ is not $c$-well-filtered.
Let $\mathcal {K}=\{\mathbb{R}\backslash C: C~\mbox{is countable}\}$. Then $\mathcal {K}\subseteq \mathcal {L}\mathcal {Q}(\mathbb{R})$. Obviously,  $\mathcal {K}$ is  countably filtered and $\bigcap\mathcal {K}=\emptyset$. However, there is no $K\in \mathcal {K}$ such that $K=\emptyset$. Therefore, $(\mathbb{R}, \tau_{coc})$ is not $c$-well-filtered.
\end{example}

Let $X$ be a $T_{0}$-space. The \emph{specialization order} $\leq$ on $X$ is defined as $x\leq y$ if and only if $x\in \mathrm{cl}({y})$. A $T_{0}$-space $X$ is called a \emph{d-space} if $X$ is a directed complete poset under the specialization order and $\mathcal {O}(X)\subseteq \sigma(X)$.

We know that each well-filtered space is a $d$-space. However, the following example shows that a $c$-well-filtered space  may not be a $d$-space.

\begin{example}\label{exam-3}
Let $L=\mathbb{N}\cup\{\omega\}$ with the order $1<2<\cdots<n<\cdots<\omega$. Then the Alexandroff space $\Gamma L=(L, a(L))$ is a $c$-well-filtered space by Proposition~\ref{prop-3.1} but not a $d$-space (see \cite[Proposition 6.5]{Xu20}).
\end{example}

Similar to well-filtered spaces, the following results hold for $c$-well-filtered spaces.

\begin{proposition}\label{prop-1}
A topological space $X$ is $c$-well-filtered if and only if for every closed subset $C$ of $X$ and each countably filtered family $\{K_{i}\}_{i\in I}$ of saturated Lindel\"{o}f subsets, if $C\cap K_{i}\neq \emptyset$ for all $i\in I$, then $\bigcap_{i\in I}K_{i}\cap C\neq\emptyset$.
\end{proposition}

\begin{proposition}\label{prop-2}
Let $(X, \tau)$ be a $c$-well-filtered $T_{0}$-space. Then $\Omega(X)$ is a countably directed complete poset and $\tau\subseteq \sigma_{c}(\Omega(X))$, where $\Omega(X)=(X, \leq_{\tau}), \leq_{\tau}$ is the specialization order of $(X, \tau)$.
\end{proposition}
\begin{proof}
Let $D$ be a countably directed subset of $\Omega(X)$. Take  $\mathcal {K}=\{{\ua_{\tau}}d: d\in D\}$, then $\mathcal {K}$ is a countably filtered family of saturated Lindel\"{o}f subsets of $X$ and $\bigcap\mathcal {K}\neq\emptyset$. Suppose that $D$ has no least upper bound. Then for any $x\in\bigcap\mathcal {K}$, there exists $y\in \bigcap\mathcal {K}$ such that $x\nleq y$. Hence, $x\in X\backslash{\da_{\tau}}y$. So we have $x\notin cl_{\tau}D$ by the fact that $(X\backslash{\da_{\tau}}y)\cap D=\emptyset$. Thus $\bigcap\mathcal {K}\cap cl_{\tau}D=\emptyset$, i.e. $\bigcap\mathcal {K}\subseteq X\backslash cl_{\tau}D$. By the  $c$-well-filteredness of $X$, we have ${\ua_{\tau}}d\subseteq X\backslash cl_{\tau}D\subseteq X\backslash D$ for some $d\in D$. This is a contradiction. Therefore, $\Omega(X)$ is a countably directed complete poset.

Let $U\in \tau$ and $D$ be a countably directed subset of $\Omega(X)$ with $\sup D\in U$. Then $\bigcap_{d\in D}\ua_{\tau}d=\ua_{\tau}\sup D\subseteq U$. Thus we have $\ua_{\tau}d\subseteq U$  for some $d\in D$ by the $c$-well-filteredness of $X$. Therefore, $U\in\sigma_{c}(\Omega(X))$.
\end{proof}

\begin{proposition}\label{prop-3}
Let $(X, \tau)$ be a $c$-well-filtered space and $A$ a saturated subset of $X$. Then the subspace $(A, \tau_{A})$ is $c$-well-filtered.
\end{proposition}
\begin{proof}
Let $\mathcal {K}$ be a countably filtered family of saturated Lindel\"{o}f subsets of $A$ and $V$ an open subset of $A$ with $\bigcap\mathcal {K}\subseteq V$. Then there exists an open subset $W$ of $X$ such that $V=W\cap A$. Hence, $\bigcap\mathcal {K}\subseteq W$. Now, we claim that $K$ is a saturated Lindel\"{o}f subset of $X$ for all $K\in \mathcal {K}$. It is easy to see that $K$ is a Lindel\"{o}f subset of $X$. So we need only to show that $K=\ua_{X}K$ for all $K\in \mathcal {K}$.
Assume $x\in \ua_{X}K$. Then there exists $y\in K$ such that $y\leq x$. As $A=\ua_{X}A$, we obtain that $x\in A$. Thus $x\in K$, which implies $K=\ua_{X}K$. So $\mathcal {K}$ is a countably filtered family of $\mathcal {L}\mathcal {Q}(X)$. By the $c$-well-filteredness of $X$, there exists $K\in \mathcal {K}$ such that $K\subseteq W$. Hence, $K\subseteq W\cap A=V$. This implies that $(A, \tau_{A})$ is $c$-well-filtered.
\end{proof}

A retract of a topological space $X$ is a topological space $Y$ such that there are two continuous mappings $r: X\rightarrow Y$ and $s: Y\rightarrow X$ such that $r\circ s=id_{Y}$.

\begin{proposition}\label{prop-4}
A retract of a $c$-well-filtered space is $c$-well-filtered.
\end{proposition}
\begin{proof}
Let $X, Y$ be topological spaces. Assume that $Y$ is $c$-well-filtered and there are continuous maps $r: Y\rightarrow X$ and $s: X\rightarrow Y$ such that $r\circ s=id_{X}$. Let $\{K_{i}\}_{i\in I}$ be a countably filtered family of saturated Lindel\"{o}f subsets of $X$ and $U$ an open subset of $X$ with $\bigcap_{i\in I}K_{i}\subseteq U$.  We know $s(K_{i})$ is a Lindel\"{o}f subset of $Y$ for all $i\in I$ by Proposition \ref{prop-2.3}. Thus $\{\ua s(K_{i})\}_{i\in I}$ is a countably filtered family of saturated Lindel\"{o}f subsets of $Y$. Next, we show that $\ua r(\ua s(K_{i}))\subseteq K_{i}$ for all $i\in I$. Take $x\in \ua r(\ua s(K_{i}))$, there exists $y\in K_{i}$ such that $r(s(y))\leq x$, i.e. $y\leq x$. So we have $x\in K_{i}$. Therefore,
$$\bigcap_{i\in I}r(\ua s(K_{i}))\subseteq \bigcap_{i\in I}\ua r(\ua s(K_{i}))\subseteq \bigcap_{i\in I}K_{i}\subseteq U.$$
Then
 $$\bigcap_{i\in I}\ua s(K_{i})\subseteq \bigcap_{i\in I} r^{-1} (r(\ua s(K_{i})))= r^{-1}(\bigcap_{i\in I} r(\ua s(K_{i}))\subseteq r^{-1}(U).$$
By the  $c$-well-filteredness of $Y$,  we have $s(K_{i})\subseteq\ua s(K_{i})\subseteq r^{-1}(U)$ for some $i\in I$. Hence,
$$K_{i}\subseteq s^{-1}(s(K_{i}))\subseteq s^{-1}( r^{-1}(U))=(r\circ s)^{-1}(U)=U.$$
So $X$ is $c$-well-filtered.
\end{proof}

A topological space $X$ is \emph{countably sober} if and only if for every countably irreducible closed subset $A$ of $X$, there exists a unique element $x\in X$ such that $A={\da}x$ (see~\cite{yang17}).

\begin{proposition}\label{prop-5}
Let $X$ be a $P$-space. If $X$ is a locally Lindel\"{o}f and $c$-well-filtered space, then $X$ is countably sober.
\end{proposition}

\begin{proof}
Let $A$ be a countably irreducible closed subset of $X$ and $\mathcal {U}=\{U\in \mathcal {O}(X): A\cap U\neq\emptyset\}$. It is easy to show that $\mathcal {U}$ is a countable filter. Let $x\in A\cap U$. Then there exists $K_{x, U}\in \mathcal {L}\mathcal {Q}(X)$ such that $x\in \mbox{int}(K_{x, U})\subseteq K_{x, U}\subseteq U$ by the fact that $X$ is locally Lindel\"{o}f.

Now, to show the proposition take $\mathcal {F}=\{K_{x, U}: x\in A\cap U\}$. Next, we show that $\mathcal {F}$ is a countably filtered family of $\mathcal {L}\mathcal {Q}(X)$.
Suppose that $\{K_{x_{i}, U_{i}}\}_{i\in \mathcal {Z}_{+}}$ is a countable family of $\mathcal {F}$. Then for any $i\in \mathcal {Z}_{+}$, $A\cap \mbox{int}(K_{x_{i}, U_{i}})\neq\emptyset$. Using the fact that $A$ is countably irreducible, we have $A\cap\bigcap_{i\in \mathcal {Z}_{+}}\mbox{int}(K_{x_{i}, U_{i}})\neq\emptyset$. Let $W=\bigcap_{i\in \mathcal {Z}_{+}}\mbox{int}(K_{x_{i}, U_{i}})$.  Since $X$ is a $P$-space, we have that $W$ is an open subset by \cite[Proposition 4.3]{yang17}. Hence, there exists $z\in A\cap W$, $K_{z, W}\in \mathcal {L}\mathcal {Q}(X)$ such that $z\in \mbox{int}(K_{z, W})\subseteq K_{z, W}\subseteq W$. So we have $K_{z, W}\in \mathcal {F}$ and $K_{z, W}\subseteq\bigcap_{i\in \mathcal {Z}_{+}}K_{x_{i}, U_{i}}$, which imply $\mathcal {F}$ is countably filtered.

Since $A\cap K_{x, U}\neq\emptyset$ for all $K_{x, U}\in \mathcal {F}$, we have $A\cap\bigcap K_{x, U}\neq\emptyset$ by Proposition \ref{prop-1}. Hence, there exists $a\in A\cap\bigcap K_{x, U}$. We only need to show that $A=\da a$. Suppose that there is a $x\in A$ such that $x\nleq a$. Thus there exists an open subset $U$ of $X$ and $x\in U$ such that $a\notin U$. Since $X$ is a locally Lindel\"{o}f space, there exists $K_{x, U}\in \mathcal {L}\mathcal {Q}(X)$ such that $x\in \mbox{int}(K_{x, U})\subseteq K_{x, U}\subseteq U$. This implies that $K_{x, U}\in \mathcal {F}$. Thus $a\in K_{x, U}$, which is a contradiction. Therefore, $A=\da a$.
\end{proof}

\section{A Hofmann-Mislove theorem}
In this section, we give a Hofmann-Mislove Theorem for $c$-well-filtered spaces.

Let $L$ be a poset and $x, y$ two elements in $L$. We say $x$ is \emph{countably way-below} $y$, written $x\ll _{c}y$, if for every countably directed subset $D$ of $L$ that has a least upper bound $\sup D$ above $y$, there is an element $d\in D$ such that $x\leq d$. Let $\da_{c}x=\{y\in L: y\ll _{c}x\}$. A countably directed complete poset $L$ is said to be a \emph{countably approximating poset} if $\da_{c}x$ is countably directed and $x=\sup\da_{c}x$ for all $x\in L$ (see~\cite{Han89}).

\begin{definition}
Let $L$ be a poset. A $\sigma$-basis $B$ of $L$ is a subset of $L$ such that for every $x\in L$, the collection $B\cap\da_{c}x$ of all elements of the $\sigma$-basis countably way-below $x$ is countably directed and $x=\sup(B\cap\da_{c}x)$.
\end{definition}
\begin{lemma}
A countably directed complete poset $L$ has a $\sigma$-basis if and only if it is countably approximating.
\end{lemma}
\begin{proof}
If $L$ is countably approximating, then $L$ is a $\sigma$-basis. Conversely, assume that $L$ has a $\sigma$-basis $B$, then $x$ is the least upper bound of $B\cap\da_{c}x$, which is a countably directed family of elements countably way-below $x$. So, it is clear that $x=\sup\da_{c}x$. Now, we claim that $\da_{c}x$ is countably directed. Let $E\in\mbox{Count}\da_{c}x$, then $e\ll _{c}x$ for every $e\in E$. Thus there exists $b_{e}\in B\cap\da_{c}x$ such that $e\leq b_{e}$ for every $e\in E$. Since $B\cap\da_{c}x$ is countably directed, there is a $b\in B\cap\da_{c}x$ such that $b_{e}\leq b$ for all $e\in E$. Hence, $E\subseteq \da b$ and $b\in\da_{c}x$. This implies that $\da_{c}x$ is countably directed and therefore we can conclude that $L$ is countably approximating.
\end{proof}
\begin{proposition}
Let $L$ be a countably approximating poset with a countable $\sigma$-basis $B$. Then $(L, \sigma_{c}(L))$ is a $c$-well-filtered space.
\end{proposition}
\begin{proof}
Suppose that $\{K_{i}\}_{i\in I}$ is a countably filtered family of saturated Lindel\"{o}f subsets of $L$ and $U$ is a $\sigma$-Scott open set with $\bigcap_{i\in I}K_{i}\subseteq U$. It is easy to show that $\{K_{i}\cap B\}_{i\in I}$ is a countably filtered family of saturated Lindel\"{o}f subsets of $B$. Hence, $\bigcap_{i\in I}(K_{i}\cap B)\subseteq U\cap B$. Then there exists $i_{o}\in I$ such that $K_{i_{0}}\cap B\subseteq U\cap B$ since $B$ with the inherited topology is $c$-well-filtered. Now, we claim that $K_{i_{0}}\subseteq U$. Assume that $K_{i_{0}}\nsubseteq U$, then there is a $t\in K_{i_{0}}$ such that $t\in L\backslash U$. Since $t=\sup(B\cap\da_{c}t)$ and $L\backslash U$ is a $\sigma$-Scott closed set, we have that $B\cap\da_{c}t\subseteq L\backslash U$. Thus, we can imply that $\ua y\cap L\backslash U \cap B\neq \O $ for every $y\in B\cap\da_{c}t$. And $\{\ua y: y\in B\cap\da_{c}t\}$ is a countably filtered family of saturated Lindel\"{o}f subsets of $B$ since $B\cap\da_{c}t$ is countably directed. So, $$\bigcap_{y\in B\cap\da_{c}t}\ua y\cap L\backslash U \cap B=\ua t \cap L\backslash U \cap B\neq \O,$$ which contradicts $t\notin L\backslash U \cap B$. Therefore, $K_{i_{0}}\subseteq U$. So, we can conclude that $(L, \sigma_{c}(L))$ is a $c$-well-filtered space.
\end{proof}
Let $L$ be a complete lattice. If $L$ is continuous, then $L$ is countably approximating. And it is clear that $\sigma(L)\subseteq \sigma_{c}(L)$. So, we can obtain that the following corollary.
\begin{corollary}
Suppose that $L$ is a continuous lattice with a countable $\sigma$-basis $B$, then $(L, \sigma(L))$ is a $c$-well-filtered space.
\end{corollary}

Let $X$ be a topological space. We discuss some properties of saturated Lindel\"{o}f subsets of $X$.

\begin{proposition}\label{prop-6}
Let $X$ be a $c$-well-filtered space. Then $K=\bigcap\mathcal {C}$ is a nonempty saturated Lindel\"{o}f set for each countable filter base $\mathcal {C}$ of nonempty saturated Lindel\"{o}f subsets of $X$. Hence, $(\mathcal {L}\mathcal {Q}(X), \supseteq)$ is a countably directed complete poset.
\end{proposition}
\begin{proof}
It is easy to see that $K\neq\emptyset$ and $K$ is saturated. Thus we only need to prove that $K$ is a Lindel\"{o}f set. Let $\mathcal {U}$ be an open cover of $K$, i.e. $K\subseteq\bigcup\mathcal {U}$. Since $X$ is $c$-well-filtered, there exists $C\in \mathcal {C}$ such that $C\subseteq\bigcup\mathcal {U}$. So there exists countably many members of $\mathcal {U}$ whose union contains $K$ by the fact that $K\subseteq C$ and $C$ is a Lindel\"{o}f subset. Therefore, $K$ is a nonempty saturated Lindel\"{o}f subset.
\end{proof}

\begin{proposition}\label{prop-7}
Let $X$ be a $P$-space.
\begin{itemize}
\item[$\mathrm{(i)}$] Let $K_{1}, K_{2}\in \mathcal {L}\mathcal {Q}(X)$ and consider the following statements:
\begin{enumerate}
\item[$(a)$] There exists $U\in \mathcal {O}(X)$ such that $K_{1}\supseteq U\supseteq K_{2}$, i.e. $\mbox{int}(K_{1})\supseteq K_{2}$;

\item[$(b)$] $K_{1}\ll _{c}K_{2}$ in $\mathcal {L}\mathcal {Q}(X)$.
\end{enumerate}

If $X$ is $c$-well-filtered, then $(a)\Rightarrow(b)$; if $X$ is locally Lindel\"{o}f, then $(b)\Rightarrow(a)$.

\item[$\mathrm{(ii)}$]  If $X$ is a locally Lindel\"{o}f and $c$-well-filtered space, then $(\mathcal {L}\mathcal {Q}(X), \supseteq)$ is a countably approximating poset.
\end{itemize}
\end{proposition}
\begin{proof}
$\mathrm{(i)}$ $(a)\Rightarrow(b)$ Suppose that $X$ is $c$-well-filtered and $\mathcal {D}$ is a countably directed subset of $\mathcal {L}\mathcal {Q}(X)$ with $K_{2}\leq \sup\mathcal {D}$. Using Proposition \ref{prop-6}, we have $\sup\mathcal {D}=\bigcap\mathcal {D}$. Hence, $K_{2}\supseteq\bigcap\mathcal {D}$. It follows from $(a)$ that there exists an open set $U$ of $X$ such that $K_{1}\supseteq U\supseteq K_{2}\supseteq\bigcap\mathcal {D}$. Thus $\bigcap\mathcal {D}\subseteq\mbox{int}(K_{1})\subseteq K_{1}$. So there exists $D\in\mathcal {D}$ such that $D\subseteq\mbox{int}(K_{1})\subseteq K_{1}$ by the $c$-well-filteredness of $X$. Therefore, $K_{1}\leq D$. This implies that $K_{1}\ll _{c}K_{2}$.

$(b)\Rightarrow(a)$ Let $U\in\mathcal {O}(X)$ and $K_{2}\subseteq U$. Since $X$ is locally Lindel\"{o}f, there exists $K_{x}\in\mathcal {L}\mathcal {Q}(X)$ such that $x\in \mbox{int}(K_{x})\subseteq K_{x}\subseteq U$ for each $x\in K_{2}\subseteq U$. Thus
$$K_{2}\subseteq \bigcup_{x\in K_{2}}\mbox{int}(K_{x})\subseteq\bigcup_{x\in K_{2}}K_{x}\subseteq U.$$
Since $K_{2}$ is a Lindel\"{o}f subset, there exists countably many members of $\{K_{x}: x\in K_{2}\}$ such that $$K_{2}\subseteq \bigcup_{i\in \mathcal {Z}_{+}}\mbox{int}(K_{x_{i}})\subseteq\bigcup_{i\in \mathcal {Z}_{+}}K_{x_{i}}\subseteq U.$$

Now, take $Q_{U}=\bigcup_{i\in \mathcal {Z}_{+}}K_{x_{i}}$ and $\mathcal {D}=\{Q_{U}: U\in\mathcal{O}(X)~\mbox{and}~ K_{2}\subseteq Q_{U}\}$. It is obvious that $\mathcal {D}\neq\emptyset$ and $Q_{U}\in \mathcal {L}\mathcal {Q}(X)$ for all $Q_{U}\in \mathcal {D}$.

We claim that $\mathcal {D}$ is countably directed. Suppose that $\{Q_{U_{i}}\}_{i\in \mathcal {Z}_{+}}$ is a countable subset of $\mathcal {D}$. Then
$$K_{2}\subseteq\bigcap_{i\in \mathcal {Z}_{+}}\mbox{int}(Q_{U_{i}})\subseteq\bigcap_{i\in \mathcal {Z}_{+}}Q_{U_{i}}\subseteq\bigcap_{i\in \mathcal {Z}_{+}}U_{i}.$$ Let $V=\bigcap_{i\in \mathcal {Z}_{+}}\mbox{int}(Q_{U_{i}})$. Then $V$ is an open subset of $X$ by~\cite[Proposition 4.3]{yang17}. From $K_{2}\subseteq V$ we know that there exists $Q_{V}\in \mathcal {L}\mathcal {Q}(X)$ such that $$K_{2}\subseteq\mbox{int}(Q_{V})\subseteq Q_{V}\subseteq V\subseteq\bigcap_{i\in \mathcal {Z}_{+}}Q_{U_{i}}.$$
This implies $Q_{V}\in\mathcal {D}$ and so $\mathcal {D}$ is countably directed.

Now we prove that $\sup\mathcal {D}=K_{2}$. It is easy to see that $K_{2}$ is an upper bound of $\mathcal {D}$. Let $K$ be another upper bound of $\mathcal {D}$. Suppose that $K_{2}\nleq K$, i.e. $K_{2}\nsupseteq K$. Then there exists $x\in K$ such that $\da x\cap K_{2}=\emptyset$. Hence, $K_{2}\subseteq X\backslash{\da}~x$. So there exists $Q\in\mathcal {L}\mathcal {Q}(X)$ such that $K_{2}\subseteq\mbox{int}(Q)\subseteq Q\subseteq X\backslash{\da}~x$, which implies $Q\in\mathcal {D}$. This contradicts $K\subseteq Q$. Therefore, $\sup\mathcal {D}=K_{2}$.

Since $K_{1}\ll _{c}K_{2}$, there exists $Q_{U}\in \mathcal {D}$ such that $K_{1}\leq Q_{U}$. Therefore, $K_{1}\supseteq Q_{U}\supseteq\mbox{int}(Q_{U})\supseteq K_{2}$.

$\mathrm{(ii)}$ Let $K\in\mathcal {L}\mathcal {Q}(X)$, $\da_{c}K=\{Q\in \mathcal {L}\mathcal {Q}(X): Q\ll _{c}K\}$. It is easy to see that $\sup\da_{c}K=K$. Hence, it is enough to show that $\da_{c}K$ is countably directed. Suppose $\{Q_{i}\}_{i\in\mathcal {Z}_{+}}$ is a countable subset of $\da_{c}K$. Then for any $i\in\mathcal {Z}_{+}$, there exists $U_{i}\in \mathcal {O}(X)$ such that $Q_{i}\supseteq U_{i}\supseteq K$ by $\mathrm{(i)}$. Thus $K\subseteq\bigcap_{i\in\mathcal {Z}_{+}}\mbox{int}(Q_{i})\subseteq \bigcap_{i\in\mathcal {Z}_{+}}Q_{i}$. Let $V=\bigcap_{i\in\mathcal {Z}_{+}}\mbox{int}(Q_{i})$.  Since $X$ is a $P$-space, we have $V$ is an open subset of $X$ and  $K\subseteq V$ by \cite[Proposition 4.3]{yang17}. So there exists $Q_{V}\in \mathcal {L}\mathcal {Q}(X)$ such that $K\subseteq\mbox{int}(Q_{V})\subseteq Q_{V}\subseteq V$. This implies that $Q_{V}\in \da_{c}K$. Therefore, we conclude that $(\mathcal {L}\mathcal {Q}(X), \supseteq)$ is a countably approximating poset by Proposition \ref{prop-6}.
\end{proof}

\begin{lemma}\label{lemm-2}
Let $X$ be a $c$-well-filtered space and $U\in \mathcal {O}(X)$. Then the set $$\phi^{'}(U)=\{K\in\mathcal {L}\mathcal {Q}(X): K\subseteq U\}$$ is a $\sigma$-Scott open countable filter in $(\mathcal {L}\mathcal {Q}(X), \supseteq)$.
\end{lemma}
\begin{proof}
It is obvious that $\phi^{'}(U)$ is an upper set. Let $\{K_{i}\}_{i\in I}$ be a countably directed subset of $\mathcal {L}\mathcal {Q}(X)$ and $\sup_{i\in I}K_{i}=\bigcap_{i\in I}K_{i}\in\phi^{'}(U)$. Then $\bigcap_{i\in I}K_{i}\subseteq U$. Hence, there exists $i_{0}\in I$ such that $K_{i_{0}}\subseteq U$ by the fact that $X$ is a $c$-well-filtered space. So $\phi^{'}(U)$ is a $\sigma$-Scott open subset. Now, we check that $\phi^{'}(U)$ is countably filtered. Let $\{K_{i}\}_{i\in \mathcal {Z}_{+}}$ be a countable subset of $\phi^{'}(U)$. Then $\inf_{i\in \mathcal {Z}_{+}}K_{i}=\bigcup_{i\in \mathcal {Z}_{+}}K_{i}$. Since $K_{i}\subseteq U$ for all $i\in \mathcal {Z}_{+}$, we have $\bigcup_{i\in \mathcal {Z}_{+}}K_{i}\subseteq U$. This implies that $\bigcup_{i\in \mathcal {Z}_{+}}K_{i}\in \phi^{'}(U)$. Therefore, $\phi^{'}(U)$ is a $\sigma$-Scott open countable filter.
\end{proof}

Let $X$ be a topological space. The set of all $\sigma$-Scott open countable filters of $(\mathcal {L}\mathcal {Q}(X), \supseteq)$ is denoted by $\mbox{OCFilt}_{\sigma}((\mathcal {L}\mathcal {Q}(X), \supseteq))$. We have the following theorem.

\begin{theorem}\label{th-2}
Let $X$ be a $P$-space. If $X$ is locally Lindel\"{o}f and $c$-well-filtered, then the mapping $$\phi^{'}: \mathcal {O}(X)\rightarrow OCFilt_{\sigma}((\mathcal {L}\mathcal {Q}(X), \supseteq)), \phi^{'}(U)=\{K\in\mathcal {L}\mathcal {Q}(X): K\subseteq U\}$$ is an order isomorphism.
\end{theorem}
\begin{proof}
$(1)$ $\phi^{'}$ is surjective. Let $\mathcal {F}\in OCFilt_{\sigma}((\mathcal {L}\mathcal {Q}(X), \supseteq))$. Take $U=\bigcup\mathcal {F}$. Now, we check that $U$ is an open subset of $X$. It is enough to show that $U$ is a neighborhood of $x$ for all $x\in U$. Assume $x\in U$. Then there exists $K\in \mathcal {F}$ such that $\ua x\subseteq K$. Using the fact that $\mathcal {F}$ is an upper set, we have $\ua x\in \mathcal {F}$. Let $\mathcal {D}=\{Q\in\mathcal {L}\mathcal {Q}(X): x\in\mbox{int}(Q)\subseteq Q\}$. Since $X$ is locally Lindel\"{o}f, we have that $\mathcal {D}\neq\emptyset$. Now, we claim that $\mathcal {D}$ is countably directed. Let $\{Q_{i}\}_{i\in \mathcal {Z}_{+}}$ be a countable subset of $\mathcal {D}$. Then $x\in \bigcap_{i\in \mathcal {Z}_{+}}\mbox{int}(Q_{i})\subseteq\bigcap_{i\in \mathcal {Z}_{+}}Q_{i}$. As $X$ is a $P$-space, we have $\bigcap_{i\in \mathcal {Z}_{+}}\mbox{int}(Q_{i})$ is an open subset of $X$ by \cite[Proposition 4.3]{yang17}. So there exists $Q\in\mathcal {L}\mathcal {Q}(X)$ such that $x\in\mbox{int}(Q)\subseteq Q\subseteq\bigcap_{i\in \mathcal {Z}_{+}}\mbox{int}(Q_{i})$, which implies that $Q\in\mathcal {D}$. Thus $\sup\mathcal {D}=\bigcap\mathcal {D}=\ua x$. Since $\mathcal {F}$ is a $\sigma$-Scott open subset, there exists $Q\in \mathcal {D}$ such that $Q\in \mathcal {F}$. Therefore, $x\in\mbox{int}(Q)\subseteq Q\subseteq U$, which implies that $U$ is a neighborhood of $x$.

Next, we show that $\phi^{'}(U)=\mathcal {F}$. Obviously, $\mathcal {F}\subseteq\phi^{'}(U)$. We only need to show that $\phi^{'}(U)\subseteq\mathcal {F}$. Let $K\in\phi^{'}(U)$ and $x\in K\subseteq U$. Then there exists $K_{x}\in \mathcal {F}$ such that $x\in \mbox{int}(K_{x})\subseteq K_{x}$. Hence, $K\subseteq\bigcup_{x\in K}K_{x}$. So there exists a countable subset $\{K_{x_{i}}: i\in \mathcal {Z}_{+}\}$ of $\{K_{x}: x\in K\}$ such that $K\subseteq\bigcup_{i\in \mathcal {Z}_{+}}K_{x_{i}}$ by the fact that $K$ is a Lindel\"{o}f subset. Thus $\bigcup_{i\in \mathcal {Z}_{+}}K_{x_{i}}\in\mathcal {F}$ because $\mathcal {F}$ is a countable filter. Therefore, $K\in \mathcal {F}$.

$(2)$ $\phi^{'}$ is an order embedding, i.e. $U_{1}\subseteq U_{2}\Leftrightarrow\phi^{'}(U_{1})\subseteq\phi^{'}(U_{2})$.

$(\Rightarrow)$ It obviously holds.

$(\Leftarrow)$ Suppose that $U_{1}\nsubseteq U_{2}$. Then there exists $u\in U_{1}$ such that $u\notin U_{2}$. Hence, $\ua u\in\phi^{'}(U_{1})$ but $\ua u\notin\phi^{'}(U_{2})$, which is a contradiction. Thus $U_{1}\subseteq U_{2}$ as desired.
\end{proof}


\bibliographystyle{./entics}

\end{document}